\documentclass{article}
\usepackage{arxiv}
\usepackage[unicode, pdftex]{hyperref}
\usepackage{url}            
\usepackage{booktabs}       
\usepackage{amsfonts}       
\usepackage{nicefrac}       
\usepackage{microtype}      
\usepackage{lipsum}		
\usepackage{listings}
\usepackage{tikz}
\usepackage{graphicx}
\usepackage{amsmath}
\usepackage{amsthm}
\usepackage{amsfonts}
\usepackage{amssymb}  
\usepackage[authoryear,round]{natbib}
\usepackage{doi}
\usepackage{color}

\title{Algorithms and topological invariants for dynamic systems. IV. Topological Equivalence of Functions and Hamiltonian Dynamical Systems}

\author{ \href{https://orcid.org/0000-0002-7164-807X}{\includegraphics[scale=0.06]{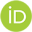}
\hspace{1mm}Alexandr O.~Prishlyak}\thanks{https://sites.google.com/view/prof-prishlyak, https://orcid.org/0000-0002-7164-807X} \\
	Department of Computer Methods of \\ 
	Mechanics and Control Processes\\
	Taras Shevchenko National University of Kyiv\\
	Kyiv, Ukraine \\
	\texttt{prishlyak@knu.ua} }

\renewcommand{\shorttitle}{Algorithms and topological invariants for dynamic systems. III. Surface Algorithms }
\newtheorem{lemma}{Lemma}
\newtheorem{theorem}{Theorem}
\newtheorem{prop}{Proposition}

\newtheorem{corol}{Corollary}
\newtheorem{definition}{Definition}

\newtheorem{remark}{Remark}

\hypersetup{
pdftitle={Algorithms and topological invariants for dynamic systems. IV. Topological Equivalence of Functions and Hamiltonian Dynamical Systems},
pdfauthor={Alexandr O. Prishlyak},
pdfkeywords={Topological classification, Morse function, Morse-Smale flow},
}

\begin{document}
\maketitle
%

\begin{abstract}
This is the fourth part of a series on algorithms and topological invariants
for dynamical systems. We study the topological equivalence of smooth
functions and of Hamiltonian flows on surfaces, and the discrete invariants,
Reeb graphs and cell decompositions, that classify them. We recall the
notions of topological equivalence and structural stability, introduce two
equivalences of functions (with and without the order of the critical values) and
show that simple Morse functions on closed orientable surfaces are classified
by their Reeb graphs with vertex odering, while simple Hamiltonian flows are classified by Reeb
graphs with marked edges. For the disc and for the sphere we obtain a recursive
description and a generating function for the numbers of such trees, and
compute the numbers of topologically non-equivalent simple Morse functions and
Hamiltonian flows with a given number of saddles. We also describe functions
with isolated, possibly degenerate, critical points by cell decompositions,
and give an algorithm and computational results for functions with three
critical points.

In the first part (arXiv:2501.15657) we discussed the basic concepts of
differential topology. In the second part (arXiv:2502.00506) we discussed the
discrete topological structures used in the topological theory of dynamical
systems. In the third part (arXiv:2502.09958) we constructed algorithms that
recognize $2$- and $3$-manifolds in simplicial complexes and regular
CW-complexes and determine the topological type of a $2$-manifold
(connectivity, orientability, genus, number of boundary components).
\end{abstract}

\keywords{Morse function \and Hamiltonian flow \and Reeb graph \and
topological equivalence \and rooted tree \and cell decomposition}

\section*{Introduction}

Functions and vector fields on manifolds are studied up to topological
equivalence by means of invariants of a discrete nature, i.e.\ objects that can
be described by a finite set of integers. For a Morse function on a surface such an
invariant is the Reeb graph, the space of connected components of the level
sets; for a Hamiltonian flow with one degree of freedom it is the same graph,
since the trajectories lie in the level sets of the Hamiltonian. Discrete
invariants are useful if there are efficient algorithms to construct,
compare and enumerate them, and the enumeration of all possible structures of
a given complexity, for example of all topological types of functions with a
given number of critical points, is itself a natural problem. Questions of
this kind, in particular the number of topologically different Morse functions
on a surface with a given number of critical points, were raised by
V.\,I.~Arnold~
and are the main subject of this part. The classification of simple Morse
functions on surfaces by Reeb graphs goes back to Kronrod and
Reeb~\cite{Kronrod1950, Reeb1946}; for further results, in particular for Hamiltonian
systems with one degree of freedom (atoms and molecules) see
the works of Fomenko, Sharko and their co-authors.

The first part~\cite{prish2025atids1} discusses the main structures used in the
topology of manifolds: vector fields, dynamical systems, Morse functions, cell
decompositions and the fundamental group. The second part~\cite{prish2025atids2}
examines discrete structures that a computer can use for specification and
manipulation. The third part~\cite{prish2025atids3} contains algorithms that
determine the topological properties of simplicial complexes and regular cell
complexes: connectivity, local planarity, orientability, the Euler
characteristic, the genus and the topological type of a surface. These
algorithms are used here to specify and recognize the cell decompositions
associated with functions.
Programs that enumerate and visualize the corresponding structures
are available as open software~\cite{prishlyak2026morsedisk, prishlyak2026handle,
prishlyak2026mssnsc, prishlyak2026rotationcodes, prishlyak2026cellembedded,
prishlyak2026chord}.

We use two equivalence relations for functions. The first, (strong)
equivalence, is the classical one: $f\circ h=h'\circ g$, where $h$ is a homeomorphism of the manifold and
$h'$ is an increasing homeomorphism of the line, so it keeps track of the order
of all critical values. The second, level equivalence, only requires a
homeomorphism that maps the components of the level sets of one function onto the
components of the level sets of the other and preserves the direction of
increase; it ignores the relative order of critical values on different
branches. Level equivalence is the one that matches the topological equivalence
of Hamiltonian flows, whose oriented trajectories are the (components of)
level sets. For simple Morse functions on closed orientable surfaces the
first equivalence is classified by Reeb graphs with the order of vertices, the
second by oriented Reeb graphs, and the equivalence of flows (where the
orientation of the surface may also change) by oriented Reeb graphs up to
reversal of all orientations. For Reeb graphs that are trees, which is the case
for the sphere and the disc, the last classification can be given in a
particularly simple form, by trees with a marked edge at every vertex of degree three.

\medskip
\noindent\textbf{Structure of the paper.}
In Section 1 
we recall the topological equivalence of
dynamical systems and diffeomorphisms, structural stability and the
Kupka--Smale theorem, and give a survey of the classifications of Morse--Smale
systems. Section 2 
is devoted to the equivalence of
functions: stability, the local classification near non-degenerate critical
points, the homotopy classification of Morse functions on surfaces, and a
complete classification, up to equivalence, of functions on the circle and on
the line with three or four extrema. In Section 3 
we define the
Reeb graph and prove its basic properties for simple Morse functions, state the
classification and realization theorems, introduce the level equivalence, and
in Subsection 3.3 
relate simple Morse functions to Hamiltonian flows
and describe the encoding of the flows by marked trees.

Section 4 
deals with Morse functions and flows on the
$2$-disc, where the Reeb graphs are rooted trees; we obtain a recurrence for the number
$K_n$ of rooted Reeb trees with $n$ saddles and compute the numbers of functions
and flows. In Section 5 
we count the Reeb graphs that are trees on
the sphere: using Otter's dissimilarity formula we express the generating
function of the numbers of simple Morse functions through $K(x)$ and obtain the numbers of
Hamiltonian flows by Burnside's lemma, with an algorithm for computing them.

In Section 6 
we consider functions with isolated critical points,
the degenerate ones included. The topological structure of such a function is
given by a cell decomposition of the surface whose $1$-skeleton is formed
by the critical levels. In Subsection 6.3 
we treat the case of three critical
points: a minimum, a maximum and a critical point whose level set is a bouquet
of $k$ circles. These functions correspond to decompositions of a closed surface
into two $k$-gons glued along a bouquet of $k$ circles, which are described by
signed permutations; we present an algorithm, a program that checks local
planarity at the vertex and determines the classes with respect to the
dihedral group $D_k$ acting on the numbering of the circles, and the resulting
numbers of topologically non-equivalent functions on orientable and
non-orientable surfaces. The last section contains the conclusions.


\section{Topological Equivalence of Dynamical Systems}\label{sec:top-equiv}

Throughout, $M$ is a closed (compact, without boundary) smooth manifold and
$r\ge 1$. We denote by $X^r(M)$ the space of $C^r$ vector fields on $M$ with
the $C^r$ topology, and by $\mathrm{Diff}^r(M)$ the space of $C^r$
diffeomorphisms with the $C^r$ topology. Since $M$ is compact, every
$X\in X^r(M)$ generates a flow $\varphi_X^t$, $t\in\mathbb R$. The
\emph{trajectory} (orbit) $\gamma_X(p)$ of $p\in M$ is the set
$\{\varphi_X^t(p)\}_{t\in\mathbb R}$, oriented by the direction of increasing
time.

\begin{definition}\label{def:te}\index{topologically equivalent vector fields}
Two vector fields $X,Y\in X^r(M)$ are called \emph{topologically}
($C^k$-smoothly, $1\le k\le r$) \emph{equivalent} if there exists a
homeomorphism ($C^k$-diffeomorphism) $h\colon M\to M$ that maps the
trajectories of $X$ onto the trajectories of $Y$ and preserves their
orientation.
\end{definition}

The map $h$ is not required to respect the time parametrization of the
trajectories; if, in addition, $h\circ\varphi_X^t=\varphi_Y^t\circ h$ for all
$t$, then $h$ is a \emph{topological conjugacy} of the flows, which is a much
stronger relation (for instance, it preserves the periods of closed
trajectories). Topological equivalence is an equivalence relation: the
inverse of $h$ and a composition of two equivalences are again equivalences.

\begin{remark}\label{rem:variants}
Definition~\ref{def:te} can be modified, and which version is meant must be
stated in every classification problem. One may additionally require that $h$
preserves the orientation of $M$ (when $M$ is oriented), or one may allow $h$
to map the oriented trajectories of $X$ either onto those of $Y$ with the
orientation preserved, or onto those of $Y$ with the orientation of
\emph{all} trajectories reversed. The latter relation means that $X$ is
equivalent to $Y$ or to $-Y$ in the sense of Definition~\ref{def:te}, and the
corresponding classification is coarser in general.
\end{remark}

\begin{theorem}\label{thm:te-props}\cite{Palis1982}
Let $h$ be a topological equivalence between the fields $X,Y\in X^r(M)$.
Then

\begin{enumerate}
\item $p\in M$ is a singular point of $X$ if and only if $h(p)$ is a singular
point of $Y$;
\item the trajectory $\gamma_X(p)$ is closed if and only if the trajectory
$\gamma_Y(h(p))$ is closed;
\item $h$ maps the $\alpha$-limit ($\omega$-limit) set of the trajectory
$\gamma_X(p)$ onto the $\alpha$-limit ($\omega$-limit) set of the trajectory
$\gamma_Y(h(p))$.
\end{enumerate}
\end{theorem}

\begin{definition}\label{def:ss}\index{structurally stable vector field}
A vector field $X\in X^r(M)$ is called \emph{structurally stable} (also
\emph{robust}; in the Russian literature, \emph{rough}) if there exists a
neighbourhood $U$ of $X$ in $X^r(M)$ such that every $Y\in U$ is
topologically equivalent to $X$.
\end{definition}

\begin{remark}
The notion of a rough system goes back to Andronov and Pontryagin, whose
original definition additionally requires the equivalence $h$ to be close to
the identity map when $Y$ is close to $X$. We use the definition above, in the
form of Peixoto and Smale. The two notions are not always the same, and we
shall say explicitly which one is used.
\end{remark}

\begin{definition}\label{def:conj}
Two diffeomorphisms $f,g\in\mathrm{Diff}^r(M)$ are called
\emph{topologically conjugate} if there exists a homeomorphism $h\colon M\to M$
such that
\[
f\circ h=h\circ g .
\]
\end{definition}

It follows that $h(g(p))=f(h(p))$ for all $p$, so $h$ maps the trajectories
$\{g^n(p)\}_{n\in\mathbb Z}$ of $g$ onto the trajectories of $f$, and a periodic point
of $g$ onto a periodic point of $f$ of the same period. Moreover
\[
h\bigl(\alpha_g(p)\bigr)=\alpha_f\bigl(h(p)\bigr),\qquad
h\bigl(\omega_g(p)\bigr)=\omega_f\bigl(h(p)\bigr),
\]
where $\omega_g(p)$ ($\alpha_g(p)$) is the set of limit points of the sequence
$g^n(p)$ as $n\to+\infty$ ($n\to-\infty$).

\begin{definition}\label{def:ss-diffeo}\index{structurally stable diffeomorphism}
A diffeomorphism $f\in\mathrm{Diff}^r(M)$ is called \emph{structurally stable}
(also \emph{robust}) if there exists a neighbourhood $U$ of $f$ in
$\mathrm{Diff}^r(M)$ such that every $g\in U$ is topologically conjugate to $f$.
\end{definition}

\medskip
Recall that a subset of a topological space is \emph{residual} if it contains
a countable intersection of open dense sets. The \emph{critical elements} of a
vector field are its singular points and closed trajectories; for a
diffeomorphism these are its periodic points. A singular point is
\emph{hyperbolic} if the linearization of the field at it has no eigenvalues
on the imaginary axis. A closed trajectory is hyperbolic if the derivative of
its Poincar\'e map has no eigenvalues of modulus~$1$, and a periodic point of
period $m$ of a diffeomorphism $f$ is hyperbolic if $D(f^m)$ at it has no
eigenvalues of modulus~$1$.

\begin{theorem}[Kupka--Smale]\label{thm:ks}\cite{Palis1982}
Let $KS$ be the set of vector fields (diffeomorphisms) that satisfy the
following conditions:
\begin{enumerate}
\item all critical elements are hyperbolic;
\item the stable and unstable manifolds of the critical elements intersect
transversally.
\end{enumerate}
Then $KS$ is residual in $X^r(M)$ (respectively, in $\mathrm{Diff}^r(M)$).
\end{theorem}

\begin{remark}
The systems of Morse--Smale type, which are defined in the following chapters,
are structurally stable. On closed orientable surfaces the converse holds
(Peixoto): structurally stable vector fields are exactly the Morse--Smale
fields, and they form an open and dense subset of $X^r(M)$. Structurally stable
diffeomorphisms of the circle $S^1$ are also exactly the Morse--Smale
diffeomorphisms, and they form an open dense subset of $\mathrm{Diff}^r(S^1)$.
In higher dimensions the situation is different: there are structurally stable
systems that are not of Morse--Smale type (e.g.\ Anosov systems), and
structurally stable systems are in general not dense in the space of all
systems; see~\cite{Palis1982}.
\end{remark}

\subsection*{Historical remarks on topological classification}

A \emph{complete topological invariant} of a class of dynamical systems is a
set of quantities (or algebraic objects) that coincide (are isomorphic) for
topologically equivalent systems of the class and differ (are not isomorphic)
for non-equivalent ones. By a \emph{topological classification} of the class we
mean the construction of a complete topological invariant, such that
equivalence (isomorphism) of invariants can be checked effectively, together
with a description of all values the invariant can take.

The topological classification of Morse--Smale vector fields on surfaces was
first obtained for the sphere by E.\,A.~Leontovich and A.\,G.~Mayer, and for
arbitrary closed surfaces by M.~Peixoto (1973), using graphs with marked
cycles. A.\,G.~Mayer and V.\,V.~Sharko used Lyapunov functions for
classification, H.~Wong used $C^*$-algebras, and G.~Fleitas and E.~Gyrik
applied circles with marked points. Morse--Smale diffeomorphisms of closed
surfaces were classified by I.~Vlasenko and, independently, by S.~Bonatti and
R.~Langevin.

On three-dimensional manifolds, G.~Fleitas studied Morse--Smale systems
without closed trajectories, and Ya.\,L.~Umanskii and O.\,O.~Prishlyak
studied those with a finite number of heteroclinic trajectories between
two-dimensional invariant manifolds. For Morse--Smale vector fields without
singular points on $S^3$, M.~Wada classified all links that can be the sets
of closed orbits of such fields. Similar questions for Smale systems (which,
unlike Morse--Smale systems, have infinitely many critical elements) were
considered by J.~Franks, M.~Sullivan and others.




\section{Equivalence of Functions}\label{sec:equiv-functions}

\subsection{Topological equivalence}

Let $M$ be a smooth manifold of dimension $n$ and let $f,g\colon M\to\mathbb R$
be smooth functions.

\begin{definition}\label{def:fe}\index{topologically equivalent functions}
The functions $f$ and $g$ are called \emph{topologically} ($C^k$-smoothly)
\emph{equivalent} if there exist homeomorphisms ($C^k$-diffeomorphisms)
$h\colon M\to M$ and $h'\colon\mathbb R\to\mathbb R$ such that
\[
f\circ h=h'\circ g .
\]
They are called \emph{strongly equivalent} if, in addition, $h'$ preserves the
orientation of the line, i.e.\ $h'$ is increasing.
\end{definition}

Both relations are equivalence relations. A homeomorphism $h'$ of $\mathbb R$ is
monotone, so in the strong equivalence $h'$ preserves the order of the values
of the functions, while in the equivalence it either preserves or reverses it;
in particular, $f$ and $-f$ are always equivalent (take $h=\mathrm{id}$,
$h'(y)=-y$), but in general not strongly equivalent.

\begin{remark}
If $M$ is oriented, one can also require that $h$ preserves its orientation.
Then there are several relations, according to whether $h$ and $h'$ preserve
or reverse the orientations of $M$ and of $\mathbb R$. Which of them is used
is stated in every classification problem; the choice corresponds to the
choice between the two kinds of equivalence of Hamiltonian flows considered in
the next chapters. The term ``conjugacy'' is not used for functions, since
it is reserved for the relation $f\circ h=h\circ g$ with the same $h$ on both
sides (Definition~\ref{def:conj}).
\end{remark}

\begin{definition}\label{def:fstable}\index{topologically stable function}
Let $M$ be closed. A function $f\in C^k(M)$ is called \emph{topologically}
(\emph{$C^k$-}) \emph{stable} if there exists a neighbourhood $U$ of $f$ in the
space $C^k(M)$ such that every $g\in U$ is topologically (respectively,
$C^k$-) equivalent to $f$.
\end{definition}

Stability is closely related to the behaviour of the critical values, because
$h'$ is a homeomorphism of the whole line: it maps the critical values of $g$
to those of $f$ bijectively and preserves their coincidences. Consequently, a
Morse function with two critical points on one level is not stable, since an
arbitrarily small perturbation separates their values. In contrast, Morse
functions whose critical values are pairwise distinct (simple Morse functions)
are topologically stable in $C^k(M)$, $k\ge2$.

Analogously to dynamical systems, the topological classification of a class of
functions is the construction of a complete topological invariant, together
with a description of all values it can take.

\subsection{Local structure}

Near a regular point every function is equivalent to a coordinate function
(the implicit function theorem), so all germs of functions at regular points
are equivalent, also in the $C^\infty$ sense.

Let $p$ be a non-degenerate critical point of $f$ of index $\lambda$. By the
Morse lemma there is a local chart in which
\[
f=f(p)-x_1^2-\dots-x_\lambda^2+x_{\lambda+1}^2+\dots+x_n^2 .
\]
It follows that two germs of functions at non-degenerate critical points $p$
and $q$, of indices $\lambda_1$ and $\lambda_2$, are locally equivalent (that
is, $f\circ h=h'\circ g$ for local homeomorphisms $h,h'$ near $p$, $q$, and
$f(p)$, $g(q)$) if and only if

\begin{enumerate}
\item $\lambda_1=\lambda_2$ in the case of strong equivalence ($h'$ increasing);
\item $\lambda_1=\lambda_2$ or $\lambda_1=n-\lambda_2$ in the case of
equivalence ($h'$ of any orientation).
\end{enumerate}

Sufficiency follows from the Morse lemma, since $-f$ has index $n-\lambda$ at the same point.
Necessity holds because the index is determined by the local homology of the
pair of sublevel sets, which is a topological invariant. On a surface
($n=2$) there are three types of non-degenerate critical points: a minimum
($\lambda=0$), a saddle ($\lambda=1$) and a maximum ($\lambda=2$). Under
strong equivalence they are pairwise non-equivalent; under equivalence a minimum
and a maximum become equivalent, while a saddle remains a saddle.

For Morse functions the homeomorphisms in a neighbourhood of a critical point can
be taken smooth, since the change of coordinates in the Morse lemma is smooth.
In general this is not so, and there are many works in singularity theory
on the possibility of reducing a function near a degenerate critical point to
a polynomial by a smooth change of coordinates.

\subsection{Global classification of functions on surfaces}

The classification of functions on the sphere by trees goes back to
G.~Reeb and A.\,S.~Kronrod. The global topological classification of Morse
functions on surfaces was obtained by V.\,V.~Sharko. A.\,T.~Fomenko introduced
atoms and molecules for the global classification of Morse functions on
surfaces.


\subsection{Homotopy of functions}

Let $M$ be a closed manifold and $k\ge2$. The space $C^\infty(M)$ with the
$C^\infty$ topology is a Fr\'echet space~\cite{Sharko1993} (for finite $k$ the space
$C^k(M)$ is a Banach space). We denote by $\mathrm{Morse}(M)$ the subspace of
Morse functions. It is open in $C^k(M)$ and locally path-connected, so its
connected components coincide with its path components.

\begin{definition}\index{Morse-homotopic functions}
Two Morse functions are called \emph{Morse-homotopic} if they belong to the same
path component of $\mathrm{Morse}(M)$, i.e.\ if one of them can be deformed
into the other continuously inside $\mathrm{Morse}(M)$.
\end{definition}

(Ordinary homotopy would be a trivial notion here, since any two functions
$M\to\mathbb R$ are homotopic.)

Let $F$ be a closed smooth surface. Denote by $\mathrm{Morse}(F,k,l)$ the
subspace of $\mathrm{Morse}(F)$ consisting of the functions with $k$ minima
and $l$ maxima.

S.\,V.~Matveev and V.\,V.~Sharko, independently and by different methods,
proved that the sets $\mathrm{Morse}(F,k,l)$ are the connected components of
$\mathrm{Morse}(F)$, thereby obtaining a homotopy classification of Morse
functions on surfaces. S.~Maksymenko generalized these results to Morse maps
of a surface to the circle: two such maps lie in the same component if and only
if they are homotopic and have the same number of critical points of each
index.


\subsection{Example: functions on the circle and on the line}\label{ex:circle-line}

We classify, up to the equivalences of Definition~\ref{def:fe}, smooth functions
with three and four critical points on the circle $S^1$ and on the line
$\mathbb R$. The critical points are not assumed to be non-degenerate, but we
assume that there are finitely many of them. On the line we consider
\emph{proper} functions, i.e.\ $|f(x)|\to\infty$ as $x\to\pm\infty$ (functions
with finite limits at infinity are discussed at the end).

\medskip\noindent\textbf{Reduction to extrema.}
Let $p$ be a critical point of $f$. Between two consecutive critical points a
continuous function without local extrema is strictly monotone, so $f$ is
strictly monotone on each arc between consecutive local extrema. If $p$ is not
a local extremum (e.g.\ $f=x^3$ near $p$), then $f$ is a homeomorphism onto an
interval near $p$, and $p$ is invisible to topological equivalence. Hence the
topological class of $f$ is determined by its local extrema only:

\begin{prop}\label{prop:extrema}
Let $M=S^1$ or $M=\mathbb R$ and let $f,g\colon M\to\mathbb R$ have finitely
many critical points (and be proper if $M=\mathbb R$). Write
$c_1,\dots,c_k$ for the local extrema of $f$ in the order along $M$, and
$d_1,\dots,d_l$ for those of $g$. Then $f$ and $g$ are strongly equivalent if and
only if $k=l$ and, for a suitable homeomorphism $h$ of $M$, the extrema
$d_i$ and $h(c_i)$ ($i=1,\dots,k$) have the same type (maximum or minimum) and
the values $g(d_1),\dots,g(d_k)$ are in the same order relation, with the same
coincidences, as the values $f(h(c_1)),\dots,f(h(c_k))$. The functions are
equivalent if this holds for $f$ and $g$ or for $f$ and $-g$.
\end{prop}

\begin{proof}
A homeomorphism preserves local extrema; if $h'$ is increasing, it preserves
their type, and if it is decreasing, it interchanges maxima and minima.
Since $h'$ is monotone, the order of the extremal values and their
coincidences are preserved (reversed if $h'$ is decreasing). Conversely, given
the stated data, choose an increasing piecewise linear homeomorphism $h'$ of $\mathbb R$ that
maps the extremal values of $g$ to those of $f$ and let $h$ be, on each arc
between consecutive extrema, the map $f^{-1}\circ h'\circ g$, where $f$ and $g$
are strictly monotone on the arc; on an unbounded arc this is again a homeomorphism
onto the corresponding arc of $f$, since $f$ and $g$ are proper. These maps
agree at the extrema and define a homeomorphism $h$ of $M$ with $f\circ h=h'\circ g$.
The case of $f$ and $-g$ is reduced to this one.
\end{proof}

\medskip\noindent\textbf{The circle.}
Maxima and minima alternate along $S^1$, so the number of local extrema is even,
and a function with exactly three critical points has two extrema and one
critical point which is not an extremum; by Proposition~\ref{prop:extrema} it is
equivalent to $\cos x$. This is the only class with two extrema, and the class with
two extrema is also the only class of functions with four critical points
two of which are not extrema.

Let $f$ have four local extrema, two maxima $M_1,M_2$ and two minima
$m_1,m_2$ in the cyclic order $M_1,m_1,M_2,m_2$. Between consecutive extrema $f$ is monotone,
so each maximum value is larger than each minimum value, and the only free
relations are those between $f(M_1)$ and $f(M_2)$ and between $f(m_1)$ and
$f(m_2)$. The rotation by two positions and the two reflections of the cyclic order that preserve the
types of the extrema allow us to interchange $M_1$ and $M_2$ and, independently,
$m_1$ and $m_2$. Therefore a function with four extrema is determined up to strong
equivalence by two yes/no data: whether $f(M_1)=f(M_2)$ and whether $f(m_1)=f(m_2)$.
We get four classes, represented by $\cos 2x$ (both equalities hold),
by $\cos 2x+\varepsilon\cos x$ (the maxima differ, the minima coincide), by
$-\cos 2x+\varepsilon\cos x$ (the minima differ, the maxima coincide), and by a
generic small perturbation of $\cos 2x$ (neither equality holds), where
$\varepsilon\ne0$ is small. Under equivalence the function $-f$ interchanges the roles of
maxima and minima, so the classes with exactly one equality coincide, and
there are \textbf{three} classes.

\medskip\noindent\textbf{The line.}
Let $f$ be proper. The extrema alternate, and the behaviour at infinity is
determined by the type of the first and the last extremum. Reflecting $x\mapsto-x$,
we may normalize the first extremum.

\emph{Three extrema.} The sequence of types is $M,m,M$ ($f\to-\infty$ at both ends)
or $m,M,m$ ($f\to+\infty$ at both ends). In each case the only free relation is
between the two extrema of the same type: whether their values coincide
(the two cases $f(M_1)<f(M_2)$ and $f(M_1)>f(M_2)$ are interchanged by the reflection).
Hence there are $2+2=4$ classes of strong equivalence, and two classes of equivalence
($f\mapsto -f$ interchanges the two sequences of types). A function with three critical points
has either three extrema (these classes) or one extremum and two critical
points that are not extrema; the latter gives the two classes of $\pm x^2$,
which coincide under equivalence.

\emph{Four extrema.} The sequence of types is $M,m,M,m$ or $m,M,m,M$, and the
reflection interchanges them, so we may assume that the first extremum is a maximum.
Denote the extrema, in the order along $\mathbb R$, by $M_1,m_1,M_2,m_2$; then
$f(M_1)>f(m_1)<f(M_2)>f(m_2)$. The free relations are those of $f(M_1)$ with
$f(M_2)$ and with $f(m_2)$, and of $f(m_1)$ with $f(m_2)$. Listing all
admissible weak orders of the four values we get \textbf{eleven} classes of strong
equivalence, given in Table~\ref{tab:line4}. The function $f\mapsto-f$ followed by
the reflection $x\mapsto -x$ preserves the normalization and acts on this set of classes as an
involution with five fixed points, so there are \textbf{eight} classes of equivalence.

\begin{table}[h]
\centering\small
\begin{tabular}{r|l|c}
No. & order of the extremal values & equivalent to No.\ under $f\mapsto -f$\\\hline
1 & $m_1=m_2<M_1=M_2$ & 1\\
2 & $m_1=m_2<M_1<M_2$ & 6\\
3 & $m_1<M_1=m_2<M_2$ & 3\\
4 & $m_1<M_1<m_2<M_2$ & 4\\
5 & $m_1=m_2<M_2<M_1$ & 8\\
6 & $m_1<m_2<M_1=M_2$ & 2\\
7 & $m_1<m_2<M_1<M_2$ & 7\\
8 & $m_2<m_1<M_1=M_2$ & 5\\
9 & $m_2<m_1<M_1<M_2$ & 10\\
10 & $m_1<m_2<M_2<M_1$ & 9\\
11 & $m_2<m_1<M_2<M_1$ & 11
\end{tabular}
\caption{Classes of proper functions on $\mathbb R$ with four local extrema
$M_1,m_1,M_2,m_2$ (in this order along the line); the rows are given up to
the normalization of the first extremum.}
\label{tab:line4}
\end{table}

A function on $\mathbb R$ with four critical points may also have fewer than
four extrema (two extrema and two other critical points, or none at all);
these functions belong to the classes with two extrema (one class) and with no
extrema (one class).

\medskip\noindent\textbf{Non-proper functions.}
If $f$ has a finite limit at an end of $\mathbb R$, this limit has to be mapped by $h'$ to the
corresponding limit of the other function, so it enters the order pattern
together with the extremal values (including its coincidence with some of them),
and the number of classes grows quickly. We do not list them here.


\section{Reeb Graph of Morse Functions on Surfaces}\label{sec:reeb}

Let $f\colon M\to\mathbb R$ be a smooth function on a smooth manifold $M$. We
define an equivalence relation on $M$: $x\sim y$ if $f(x)=f(y)$ and $x$, $y$
lie in the same connected component of the level set $f^{-1}(f(x))$.

\begin{definition}\label{def:reeb}\index{Reeb graph}
The \emph{Reeb graph} of $f$ is the quotient space
\[
\Gamma_f=M/{\sim}
\]
with the quotient topology. Let $\pi\colon M\to\Gamma_f$ be the projection. Since $f$ is
constant on the classes of $\sim$, it induces a continuous function
$\tilde f\colon\Gamma_f\to\mathbb R$ with $f=\tilde f\circ\pi$.
\end{definition}

Thus the points of the Reeb graph are the connected components of the level
sets of $f$. For the functions considered below $\Gamma_f$ is a finite graph:
its vertices are the images of the connected components of the critical
levels that contain critical points, and its edges are the components of the complement
of the vertices. The function $\tilde f$ is strictly monotone on each edge,
and the direction of increase of $\tilde f$ orients every edge.

\begin{definition}\label{def:simple}\index{simple Morse function}
A function $f$ on a closed manifold is a \emph{Morse function} if all its
critical points are non-degenerate. A Morse function is called \emph{simple} if
its critical values are pairwise distinct, i.e.\ every critical level contains
exactly one critical point.
\end{definition}

(By the previous section, among Morse functions the topologically stable ones are exactly the simple ones.)
In this section all surfaces are connected. For a simple Morse function the
values of $\tilde f$ at the vertices are pairwise distinct, so
$v_i<v_j$ if $\tilde f(v_i)<\tilde f(v_j)$ is a \emph{strict total order} on the
set of vertices. Every edge is oriented from its smaller endpoint to the larger
one.

\begin{lemma}\label{lem:structure}
Let $f$ be a simple Morse function on a closed surface $M$.
\begin{enumerate}
\item $\Gamma_f$ is a finite connected graph, and each open edge has as its
preimage in $M$ an open cylinder $S^1\times(0,1)$ (on which $f$ is the
coordinate on the interval).
\item A local extremum of $f$ is a vertex of valence $1$: a source if it is a
minimum and a sink if it is a maximum.
\item Let $v$ be the vertex of a saddle. If a neighbourhood of the
singular component is orientable (in particular, if $M$ is orientable), then $v$ has
valence~$3$ with one incoming and two outgoing edges, or two incoming and one
outgoing edge; in particular, $v$ is neither a source nor a sink.
If the neighbourhood is non-orientable, then $v$ has valence~$2$, with one incoming
and one outgoing edge.
\item $\Gamma_f$ contains no directed cycles.
\end{enumerate}
\end{lemma}

\begin{proof}
(1) Over an open edge the function $f$ has no critical points, so the
gradient flow trivializes the preimage as a product of a regular level
component (a circle) and an interval. (2) A non-degenerate extremum has a
disc neighbourhood bounded by a level circle. (3) Let $c$ be the critical value and
let $N$ be the component of $f^{-1}[c-\varepsilon,c+\varepsilon]$ containing the saddle.
It retracts onto a figure eight, so $\chi(N)=-1$. Its boundary consists of
circles of $\{f=c-\varepsilon\}$ and of $\{f=c+\varepsilon\}$, and both kinds are present,
because near a saddle $f$ takes values both below and above $c$. Hence $N$ has at
least two boundary circles. An orientable surface with $\chi=-1$ and at
least two boundary circles is a pair of pants (three boundary circles: one
below and two above, or two below and one above); a non-orientable one is a
Möbius band with a hole (two boundary circles, one below and one above). (4) $\tilde f$
increases along the edges.
\end{proof}

\begin{corol}\label{cor:genus}
Let $M$ be a closed orientable surface of genus $g$ and $f$ a simple Morse
function with $s$ saddles, $L=s+2-2g$ extrema. Then $\Gamma_f$ has $V=2(s+1-g)$
vertices and $E=2s+1-g$ edges, and the first Betti number of $\Gamma_f$ is $g$.
In particular $V$ is even.
\end{corol}

\begin{proof}
By Morse theory, $L-s=\chi(M)=2-2g$. The
degrees of the vertices are $3$ at the saddles and $1$ at the extrema, so
$2E=3s+L$, and $b_1(\Gamma_f)=E-V+1=g$.
\end{proof}


\begin{theorem}\label{thm:classification}\cite{Kronrod1950, Reeb1946}
Let $f,g$ be simple Morse functions on a closed orientable surface.
\begin{enumerate}
\item $f$ and $g$ are strongly equivalent (Definition~\ref{def:fe}) if and only if there
is an isomorphism $\Gamma_f\to\Gamma_g$ of graphs that preserves the order on
the sets of vertices.
\item $f$ and $g$ are equivalent if and only if there is an isomorphism
$\Gamma_f\to\Gamma_g$ of graphs that either preserves or reverses the order
on the sets of vertices (in the second case it reverses the orientations of all edges).
\end{enumerate}
\end{theorem}


There is also a realization theorem.

\begin{theorem}\label{thm:realization}\cite{Kronrod1950, Reeb1946}
Let $G$ be a finite connected oriented graph such that
\begin{enumerate}
\item all vertices have valence $1$ or $3$;
\item the vertices of valence $3$ are neither sinks nor sources;
\item $G$ contains no directed cycles.
\end{enumerate}
Then every strict total order on the vertices of $G$ compatible with the
orientation of the edges is realized by a simple Morse function on a closed
orientable surface of genus $b_1(G)$, that is, there is a simple Morse
function $f$ and an isomorphism $\Gamma_f\to G$ that preserves the order.
\end{theorem}

Such a total order always exists, since an oriented graph without directed
cycles has a topological ordering; hence every graph satisfying 1--3 is the
Reeb graph of some simple Morse function.

\subsection{Level equivalence}

The order on the vertices in Theorem~\ref{thm:classification} records the
relative position of the critical values of vertices that are not
connected by a path in $\Gamma_f$ monotone along the orientation. The
equivalence of Definition~\ref{def:fe} sees this order because $h'$ is a
homeomorphism of the whole line. For Hamiltonian flows, whose trajectories
are the components of the level sets, this information is irrelevant, so we also
consider the following coarser relation.

\begin{definition}\label{def:level}\index{level equivalence}
Two functions $f,g$ on a surface $M$ are called \emph{level equivalent} if
there is a homeomorphism $h\colon M\to M$ that maps the connected components of
the level sets of $g$ onto the connected components of the level sets of $f$
and preserves the direction of increase of the functions, i.e.\ induces a
homeomorphism $\Gamma_g\to\Gamma_f$ that maps oriented edges onto oriented edges.
\end{definition}

\begin{theorem}\label{thm:level}
Two simple Morse functions on a closed orientable surface are level equivalent
if and only if their Reeb graphs are isomorphic as \emph{oriented} graphs,
without any condition on the order of the vertices.
\end{theorem}

\begin{proof}
A level equivalence induces an isomorphism of oriented graphs. Conversely,
let $\varphi\colon\Gamma_g\to\Gamma_f$ be an isomorphism of oriented
graphs. The order of vertices induced from $f$ by $\varphi$ is a strict total
order on the vertices of $\Gamma_g$ compatible with the orientation. Choose a
continuous function $u\colon\Gamma_g\to\mathbb R$, strictly increasing along the edges,
with the values at the vertices realizing this order. The function
$g'=u\circ\pi_g$ has the same level components as $g$ (they are the fibres of
$\pi_g$), and by Theorem~\ref{thm:classification} (applied to $g'$ and $f$, whose Reeb graphs are
now isomorphic as ordered graphs) there is a homeomorphism $h$ and a
homeomorphism $h'$ of $\mathbb R$ with $f\circ h=h'\circ g'$. This $h$ maps the level components of $g$ onto those of $f$ and preserves the
direction of increase.
\end{proof}

\begin{remark}
Level equivalence is strictly weaker than strong equivalence: two simple Morse
functions with isomorphic oriented Reeb graphs need not be strongly
equivalent. The smallest example occurs on the sphere with three saddles:
take the tree with a middle saddle $s_2$ having two incoming edges, from the saddles
$s_1$ and $s_3$, and one outgoing edge to a maximum; let $s_1$ have a minimum and a maximum as the remaining
neighbours, and $s_3$ two minima. The saddles $s_1$ and $s_3$ are not
comparable along the orientation, and the orders $s_1<s_3$ and $s_3<s_1$ give
two strongly non-equivalent functions (the tree has no automorphism
exchanging $s_1$ and $s_3$). In the following sections the classes of simple
Morse functions on the sphere are counted up to level equivalence, and the
classes of Hamiltonian flows up to level equivalence and the simultaneous reversal of
all orientations.
\end{remark}

\subsection{Non-orientable surfaces and surfaces with boundary}

For closed non-orientable surfaces the Reeb graph does not contain complete
information about the function: there exist topologically inequivalent
functions with isomorphic Reeb graphs. By Lemma~\ref{lem:structure}, the Reeb
graph on a non-orientable surface may have vertices of valence $2$ (one
incoming and one outgoing edge). If the graph has such a vertex, the surface
is non-orientable, and the graph determines its topological type, because the
Euler characteristic is the number of extrema minus the number of saddles. If
there are no such vertices and the graph is not a tree, the graph does not
determine whether the surface is orientable: for example, on the torus and on
the Klein bottle there are simple Morse functions with four critical points
and isomorphic Reeb graphs (a minimum, two saddles joined by a double edge, a
maximum). For a tree without vertices of valence $2$ the surface is the sphere.

For compact surfaces with boundary the components of the preimages of regular
values may be not only circles but also segments (with endpoints on the boundary),
so the corresponding edges of the Reeb graph can be coloured in two colours.




\subsection{Hamiltonian flows and Morse functions on surfaces}\label{subsec:hamiltonian}

Let $M$ be an oriented smooth surface with an area form $\omega$ compatible with
the orientation (for example, the area form $\omega_g$ of a Riemannian metric
$g$). For a smooth function $H\colon M\to\mathbb R$ the \emph{Hamiltonian
vector field} $X_H$ is defined by
\[
\iota_{X_H}\omega=dH ,
\]
where $\iota_{X_H}$ is the contraction with $X_H$. (With the opposite sign
convention, $\iota_{X_H}\omega=-dH$, all trajectories are reversed, which
amounts to replacing $H$ by $-H$.)

\begin{lemma}\label{lem:hamiltonian}
Let $X_H$ be the Hamiltonian vector field of $H$.
\begin{enumerate}
\item $dH(X_H)=\omega(X_H,X_H)=0$, so every trajectory lies in a level set of
$H$. The singular points of $X_H$ are exactly the critical points of $H$.
\item If $c$ is a regular value, then every connected component of
$H^{-1}(c)$ is a closed trajectory. Every connected component of a critical
level is a union of trajectories: singular points and separatrices.
\item If $\omega'=\varphi\omega$ with $\varphi>0$ is another area form
compatible with the orientation, then $X'_H=X_H/\varphi$. Hence the oriented
trajectories do not depend on the choice of the area form, in particular on the
Riemannian metric.
\item Let $p$ be a regular point and let $v$ be a vector with $dH(v)>0$. Then
$(X_H,v)$ is a positively oriented basis. Hence the direction of the trajectory
through $p$ is determined by the orientation of $M$ and the direction of
increase of $H$.
\end{enumerate}
\end{lemma}

\begin{proof}
(1) The form $\omega$ is non-degenerate, so $X_H(p)=0$ if and only if $dH(p)=0$.
(2) A compact connected one-dimensional level component without singular
points is a circle that is a single trajectory. (3) $\iota_{X'}\omega'=
\varphi\,\iota_{X'}\omega=dH$ gives $X'=X/\varphi$. (4) $\omega(X_H,v)=dH(v)>0$.
\end{proof}

Thus a Morse function on a closed oriented surface determines a Hamiltonian flow
whose trajectories lie in the connected components of the level sets: the regular
components are closed trajectories, and each component of a critical level is
formed by a saddle (or an extremum) and its separatrices.

\begin{definition}\label{def:simple-flow}\index{simple Hamiltonian flow}
A Hamiltonian flow $X_H$ of a Morse function $H$ on a closed surface is called
\emph{simple} if every connected component of every critical level of $H$ contains exactly
one critical point.
\end{definition}

This is a property of the oriented trajectories only: no two saddles are joined
by a separatrix. It does not require the critical values to be distinct. For
such a flow the Reeb graph of $H$ satisfies the conditions of
Theorem~\ref{thm:realization}, so the flow is equivalent to the Hamiltonian flow
of a simple Morse function with the same Reeb graph (Theorems~\ref{thm:realization}
and~\ref{thm:level}); we can therefore assume that $H$ is simple.

\begin{prop}\label{prop:flows}
Let $H_1,H_2$ be simple Morse functions on closed oriented surfaces $M_1,M_2$,
and let $X_1,X_2$ be their Hamiltonian flows (for any area forms).
\begin{enumerate}
\item $X_1$ and $X_2$ are topologically equivalent by an orientation-preserving
homeomorphism $M_1\to M_2$ if and only if the Reeb graphs $\Gamma_{H_1}$ and
$\Gamma_{H_2}$ are isomorphic as oriented graphs. Equivalently, $H_1$ and $H_2$ are level
equivalent (Definition~\ref{def:level}) by an orientation-preserving homeomorphism.
\item $X_1$ and $X_2$ are topologically equivalent (by an arbitrary homeomorphism)
if and only if $\Gamma_{H_1}$ is isomorphic to $\Gamma_{H_2}$ or to
$\Gamma_{H_2}$ with all orientations reversed.
\end{enumerate}
\end{prop}

The proof uses the following lemma, which shows that the orientation of the
surface carries no additional information.

\begin{lemma}\label{lem:reflection}
Let $H$ be a simple Morse function on a closed orientable surface $M$. There is an
orientation-reversing homeomorphism $\rho$ of $M$ with $H\circ\rho=H$.
\end{lemma}

\begin{proof}[Sketch]
By Lemma~\ref{lem:structure}, $M$ is glued from standard pieces along level
circles: discs at the extrema, cylinders over the edges and pairs of pants at
the saddles. Each piece has an involution that preserves $H$, reverses the
orientation and acts on every boundary circle as a reflection. The gluing
homeomorphisms of the boundary circles can be chosen to commute with these
reflections. By Theorem~\ref{thm:classification} the glued surface with the
glued function is equivalent to $(M,H)$, and the involutions combine to a
required $\rho$.
\end{proof}

\begin{proof}[Proof of Proposition~\ref{prop:flows}]
Let $h\colon M_1\to M_2$ be a homeomorphism that maps oriented trajectories of
$X_1$ onto oriented trajectories of $X_2$. Then $h$ maps closed trajectories
to closed trajectories and critical level components to critical level components
(by Theorem~\ref{thm:te-props}), hence level components of $H_1$ onto level components of
$H_2$. By Lemma~\ref{lem:hamiltonian}(4) the direction of a trajectory is given by
the pair (orientation of the surface, direction of increase of $H$), and $h$
preserves this direction. Therefore $h$ either preserves both the orientation of the surface and the direction
of increase, or reverses both. In the first case $h$ induces an isomorphism of
oriented Reeb graphs, in the second case an isomorphism of $\Gamma_{H_1}$ onto $\Gamma_{H_2}$
with all orientations reversed. This proves the necessity in both statements.

Conversely, let $\Gamma_{H_1}\cong\Gamma_{H_2}$ be an isomorphism of oriented
graphs. By Theorem~\ref{thm:level} there is a level equivalence $h$ of $H_1$ and $H_2$, and by
Lemma~\ref{lem:reflection} we may compose $h$ with $\rho$ and assume that $h$ preserves the orientation. By
Lemma~\ref{lem:hamiltonian}(4) it then maps oriented trajectories onto oriented
trajectories (on a critical component $h$ maps the saddle to the saddle and the
separatrices to separatrices, since the topological type of the point is preserved). If
$\Gamma_{H_1}$ is isomorphic to $\Gamma_{H_2}$ with the orientations reversed, apply the same
argument to $H_1$ and $-H_2$ and then compose with $\rho$ in order to get an orientation-reversing map that
reverses the direction of increase; it again preserves the direction of
the trajectories. 
\end{proof}

Thus, if the equivalences are required to preserve the orientation of the surface (and
the direction of trajectories, in the case of flows), then the classification
of simple Morse functions up to level equivalence is the same as the
classification of the corresponding Hamiltonian flows: both are given by oriented
Reeb graphs up to isomorphism. If the orientation of the surface is allowed to
change, the Hamiltonian flows are classified by oriented Reeb graphs up to
isomorphism and the simultaneous reversal of all orientations.

\subsection*{Encoding by marked edges}

For a general Morse function the orientation of the edges of the Reeb graph is
induced by the direction of increase of the function. For simple Hamiltonian flows
it is more convenient to use a different encoding, which does not depend on the
choice of the direction. At each vertex of degree~3 (a saddle) exactly one of the
three incident edges is directed differently from the other two (it is the only
incoming edge or the only outgoing edge). We \emph{mark} this edge at the vertex. A
\emph{marked tree} is a tree with vertices of degree~$1$ and $3$ together with a
choice of one marked edge at each vertex of degree~$3$ (an edge may be marked at one of its
ends and not at the other). Two marked trees are isomorphic if there is an
isomorphism of trees that maps marked edges to marked edges.

\begin{prop}\label{prop:marked}
Let $T$ be a tree with vertices of degree $1$ and $3$. For every
marking of $T$ there are exactly two orientations of the edges of $T$ such
that at every saddle the marked edge is directed oppositely to the other two incident
edges (opposite with respect to the vertex: one is incoming and the other two are outgoing, or vice versa);
these two orientations are reverse to each other. Consequently, isomorphism classes of marked trees with
$n$ saddles correspond bijectively to isomorphism classes of oriented trees
satisfying the conditions of Theorem~\ref{thm:realization} up to isomorphism and
the simultaneous reversal of all orientations.
\end{prop}

\begin{proof}
Fix an edge $e_0$ (for instance, the edge incident to the root) and one of its two orientations.
The orientation of the other edges is then determined step by step along the
tree. If an edge $e$ at a vertex $v$ has a known direction relative to $v$
(incoming or outgoing), then the directions of the two other edges at $v$ are
determined by the rule: if $e$ is the marked edge, the two others are
directed oppositely to it; if $e$ is not marked, then the other unmarked edge has the same
direction as $e$, and the marked edge has the opposite direction. The direction relative to the
other end of an edge is the opposite one. Since $T$ has no cycles, the
choices never conflict, and there are exactly two orientations (that
differ by the choice at $e_0$), reverse to each other. At each saddle the result
is neither a sink nor a source. Conversely, from an oriented tree the marking is
recovered by marking at each saddle the edge directed differently from the other two. The
maps are inverse to each other, and an isomorphism of marked trees carries the
orientations to each other or to the reversed ones, which gives the
bijection of the classes.
\end{proof}

\begin{remark}
By Proposition~\ref{prop:flows}, for trees (the Reeb graphs of functions on the
sphere and on the disc) the classes of simple Hamiltonian flows up to equivalence correspond to
marked trees up to isomorphism, and the classes of simple Morse functions up to level equivalence correspond
to oriented trees up to isomorphism. The two counting problems are therefore the counting
of marked trees and of oriented trees with $n$ saddles.
\end{remark}

\section{Simple Morse Functions and Hamiltonian Flows on the 2-Disk}

To determine the number of topologically non-equivalent simple Hamiltonian flows, we count the Reeb graphs with a given number of vertices corresponding to saddle critical points.

We begin with the case of trees having vertices of degrees one and three, with one of the degree-one vertices distinguished as the root. Such trees correspond to Reeb graphs of Morse functions on the 2-disk, where the root represents the boundary fiber.

At each vertex of degree three, we mark one of the three incident edges. As explained in the previous section, this marking determines the orientation pattern of the three incident edges: the marked edge has the opposite orientation to the other two.

We orient the edges of the tree according to the following rules:
\begin{enumerate}
\item The edge incident to the root is directed away from the root.
\item If the marked edge at a saddle vertex is directed towards that vertex, the other two edges are directed away from it.
\item If the marked edge at a saddle vertex is directed away from that vertex, the other two edges are directed towards it.
\end{enumerate}

Since the graph is a tree, these rules uniquely determine the orientation of all its edges. Indeed, starting with the edge incident to the root, we can successively determine the orientations of all edges adjacent to the already oriented ones. The process terminates when every edge of the tree has been oriented.

\begin{lemma}
The number $K_k$ of rooted Reeb trees with $k$ saddle vertices (vertices of type $A_0$), all other vertices being of type $C$, can be computed using the following recursive formulas:
\[
K_0=1,
\]
for even $k=2n$,
\begin{equation}
K_{2n}=3\sum_{i=0}^{n-1}K_iK_{2n-1-i},
\end{equation}
and for odd $k=2n+1$,
\begin{equation}
K_{2n+1}=
\frac{3K_n^2+K_n}{2}
+3\sum_{i=0}^{n-1}K_iK_{2n-i}.
\end{equation}
\end{lemma}

\begin{proof}
If k=0, then the tree has no saddles. Such a graph is unique. It has two vertices and one edge connecting them. From now on, we  assume that \( k > 0 \).

Consider the saddle $v$ closest to the root of a rooted tree $G$ with $k$ saddle points (vertices of type $A_0$) and other vertices of type $C$. It is a vertex of degree 3. We cut the graph by $v$. We obtain three connected components. One of them is the edge $e$ connecting $v$ to the root. The other two components ($A$ and $B$) can be viewed as rooted Reeb graphs, where the root is $v$. Let us assume that the number of vertices $|V_A|$ of  no fewer than the number of vertices $|V_B|$ of $B$ ($|V_A|\geq |V_B|$).

Then there are three possibilities:

1) the marked edge to a vertex $v$ is $e$; 

2) the marked edge belongs to  $A$; 

3) the marked edge belongs to $B$.




Multiplying by 3 in formulas (2) and (3) is to indicate which subgraph with the marked edge is connected to.

If $k=2n$, then the number of graph vertices is different ($|V_A|>|V_B|$), and $A$ cannot be isomorphic to $B$. The number of non-isomorphic trees with $i$ saddles is equal to $K_i$, and with $2n - 1 - i$ saddles is equal to $K_{2n - 1 - i}$, so the total number of trees $G$ of this type is equal to $K_i K_{2n - 1 - i}$. By considering all possible values of $i$, we obtain formula (2). 

If $k = 2n + 1$ and $|V_A|>|V_B|$, then we have similar products as in the case of $k = 2n$. Let's consider the case of $|V_A|=|V_B|=n$ in more detail. If edge $e$ is not marked, then we relabel the graphs so that the marked edge belongs to graph $A$. Since each of the graphs can be chosen in $K_n$ ways, there are $K_n^2$ graphs $G$ of this type. If edge $e$ is marked, then two cases are possible: 1) graph $A$ is not isomorphic to graph $B$, 2) graphs $A$ and $B$ are isomorphic. In the first case, we have $K_n(K_n-1)/2$ non-isomorphic graphs, and in the second case, $K_n$ graphs. Summing everything together, we obtain formula (3).
\end{proof}
The results of the calculations of the number of Reeb rooted trees on $\mathbb{R}P^2$ are presented in Table \ref{PP}.





The results of the calculations of the number of simple Hamiltonian flows on $D^2$ are presented in Table \ref{PP}.

$ $

	\begin{table}[ht] \footnotesize
		\caption{\textbf{Numbers of rooted trees and simple Hamiltonian flows on $D^2$}}\vspace{0pt}
		\centering
		\begin{tabular}{|c|c|c|c|c|c|c|c|c|c|} 
						\hline
  $k$ &0& 1& 2&3 &4 &5 &6 &7 & 8 \\
\hline
$K_k$ & 1& 2& 6 & 25 & 111&540 &2736 &14396 &77649 \\
\hline
		\end{tabular}
		\label{PP}
	\end{table}

\begin{remark}
Morse flows on the disc with holes, including those whose Reeb graphs are the
trees counted above, can be encoded by codes, enumerated for given numbers of
vertices and edges, and visualized with~\cite{prishlyak2026morsedisk}.
\end{remark}

\section{Reeb trees on the sphere}


By a \emph{Reeb tree} we mean a finite tree $T$ whose edges are oriented
and whose vertices have degree $1$ or $3$, such that every vertex of
degree $3$ (a \emph{saddle}) has one or two incoming edges (and hence two or
one outgoing edges, respectively). The vertices of degree $1$ are
\emph{leaves}: a leaf is a \emph{source} (a minimum) if its edge is outgoing,
and a \emph{sink} (a maximum) if its edge is incoming.
If $T$ has $n$ saddles, then $T$ has $L$ leaves, $V=n+L$ vertices and $V-1$ edges,
and counting degrees gives $3n+L=2(V-1)$, i.e.
\[
L=n+2,\qquad V=2n+2,\qquad E=2n+1 .
\]
Every orientation of a tree is acyclic, so it is induced by a function
on the vertex set that increases along every edge. Hence every Reeb tree is the Reeb graph of
a simple Morse function on $S^2$ (a function with distinct critical values is not
required, since the Reeb graph does not depend on the order of the
critical values on non-adjacent branches).

An \emph{isomorphism} of Reeb trees $T\to T'$ is a graph isomorphism
that maps every edge to an edge with the same orientation.
A \emph{reversing isomorphism} is a graph isomorphism that reverses the orientation
of every edge. Let $\overline T$ denote the tree $T$ with all
orientations reversed (reversal turns sources into sinks and
one-incoming saddles into two-incoming saddles, so $\overline T$ is again a Reeb tree).
Clearly $T'\cong \overline T$ if and only if there is a reversing isomorphism $T\to T'$.
A reversing automorphism of $T$ is called an \emph{anti-automorphism}.

Let $T_n$ be the number of isomorphism classes of Reeb trees with $n$ saddles.
These are the Reeb graphs of simple Morse functions on $S^2$ with $n$ saddles up to
equivalence. Let $H_n$ be the number of classes of Reeb trees with $n$ saddles
up to isomorphisms and reversing isomorphisms;
these are the Reeb graphs of simple Hamiltonian flows on $S^2$ with
$n$ saddles. Reversing all orientations is an involution on the set of $T_n$ classes,
and $H_n$ is the number of its orbits.

\medskip\noindent\textbf{Branches.}
A \emph{branch} is a Reeb tree $B$ with a marked leaf $r$ (its \emph{root})
together with the \emph{connecting edge} $e_B$ incident to $r$,
considered up to isomorphisms preserving $r$.
(Equivalently, it is a Reeb tree with one half-open edge attached to
the root; a branch with $k$ saddles is obtained from a ``rooted Reeb tree'' of the
previous section by choosing the type of the root.)
A branch is of type $+$ if $e_B$ is directed away from the attaching point,
that is, the root is a sink of the branch, and of type $-$ otherwise.
In other words, if a branch hangs on an edge $u\to v$,
then the branch at the end $v$ has type $+$ and the branch at the end $u$ has type $-$.
Reversal of orientations is a bijection between branches of type $+$ and branches of type $-$
with the same number of saddles. By Lemma~1 their number is $K_k$ for each type
(for the type of the root fixed). Put
\[
K(x)=\sum_{k\ge0}K_kx^k .
\]

\medskip\noindent\textbf{Symmetric square.}
For a power series $F(x)$ with non-negative integer coefficients that counts some
objects, the series counting \emph{unordered} pairs of such objects
(equal objects allowed) is $\tfrac12\bigl(F(x)^2+F(x^2)\bigr)$.

\textbf{The recursion in generating-function form}

\begin{lemma}\label{lem:gf}
The series $K(x)$ satisfies
\begin{equation}\label{eq:K}
K(x)=1+\frac{x}{2}\Bigl(3K(x)^2+K(x^2)\Bigr).
\end{equation}
This is equivalent to the recursion of Lemma~1.
\end{lemma}

\begin{proof}
Let $B$ be a branch of type $+$ with $k\ge1$ saddles. Its root is a leaf and
the neighbouring vertex $v$ is a saddle, whose incoming edge is the connecting edge $e_B$.
The other two edges of $v$ lead to two sub-branches, which form an unordered pair.
The saddle $v$ has one or two incoming edges, so at most one of the
two other edges is incoming. A sub-branch has type $+$ if its edge is outgoing from $v$
and type $-$ if it is incoming. Hence the allowed pairs of types are $\{+,+\}$ and
$\{+,-\}$, and the pair $\{-,-\}$ is forbidden. The number of branches with $k-1$ saddles in total equals
\[
[x^{k-1}]\Bigl(\tfrac12\bigl(K(x)^2+K(x^2)\bigr)+K(x)^2\Bigr)
=[x^{k-1}]\tfrac12\bigl(3K(x)^2+K(x^2)\bigr),
\]
where the first term counts pairs $\{+,+\}$ (symmetric square) and the second term counts
pairs $\{+,-\}$ (types distinguish the two elements, so the pair is ordered).
For $k=0$ the branch is a single leaf, which gives the term $1$.
Branches of type $-$ are treated in the same way (the forbidden pair is $\{+,+\}$),
which gives the same series. Separating the coefficients of even and odd degree
in \eqref{eq:K}, we get the two formulas of Lemma~1.
\end{proof}


\begin{theorem}\label{thm:main}
\begin{align}
\sum_{n\ge 0}T_nx^n &= K(x)^2\bigl(1-2xK(x)\bigr),\label{eq:T}\\
H_n &= \frac12\bigl(T_n+F_n\bigr),\qquad
F_n=\begin{cases}K_{n/2},& n \text{ even},\\ 0,& n\text{ odd}.\end{cases}\label{eq:H}
\end{align}
Equivalently,
\[
T_n=\sum_{i+j=n}K_iK_j-2\sum_{i+j+l=n-1}K_iK_jK_l .
\]
\end{theorem}

\begin{proof}
\emph{The number $T_n$.}
We use Otter's dissimilarity theorem for trees: the number of isomorphism classes of trees
equals the number of classes of vertex-rooted trees minus the number of classes of
edge-rooted trees plus the number of classes of trees with an automorphism that
interchanges the two ends of an edge. In our setting all isomorphisms preserve edge orientations,
so no automorphism interchanges the ends of an edge, and the last term is zero.
Thus $T_n=R_n-E_n$, where $R_n$ and $E_n$ count the classes of vertex-rooted and edge-rooted
Reeb trees with $n$ saddles.

\emph{Edge-rooted trees.} Let $e=(u\to v)$ be the marked edge.
Cutting $e$ gives a branch at $u$ (of type $-$) and a branch at $v$ (of type $+$).
An isomorphism preserves $e$ and its direction, so the two branches are determined independently,
and the generating function is
\[
\sum E_nx^n=K(x)^2 .
\]

\emph{Vertex-rooted trees.} If the root is a leaf, then it is a source or a sink,
and the rest of the tree is one branch, which gives $2K(x)$.
If the root is a saddle, then it carries an unordered triple of branches.
The saddle has one or two outgoing edges, so the triple contains
one or two branches of type $+$ (an outgoing edge gives a branch of type $+$).
The triples of types $\{+,+,-\}$ and $\{+,-,-\}$ give
$K(x)\cdot\tfrac12(K(x)^2+K(x^2))$ each, and
\[
\sum R_nx^n=2K+xK\bigl(K^2+K(x^2)\bigr).
\]
By Lemma~\ref{lem:gf}, $xK(x^2)=2(K-1)-3xK^2$, so
\[
\sum R_nx^n=2K+xK^3+K\bigl(2K-2-3xK^2\bigr)=2K^2-2xK^3 .
\]
Hence $\sum T_nx^n=2K^2-2xK^3-K^2=K^2(1-2xK)$, which is \eqref{eq:T}.

\emph{The number $H_n$.}
By Burnside's lemma applied to the reversal involution,
$H_n=\tfrac12(T_n+F_n)$, where $F_n$ is the number of classes with
$T\cong\overline T$, i.e.\ of Reeb trees having an anti-automorphism $\sigma$.
The map $\sigma$ is an automorphism of the underlying tree, so it fixes
the centre of $T$ (a vertex or an edge). We consider the cases.
\begin{itemize}
\item A leaf cannot be fixed: it is the end of one edge, which would be mapped to itself with
both ends fixed, but $\sigma$ reverses its orientation.
\item A saddle $v$ cannot be fixed: $\sigma$ permutes its three edges and maps
$k$ incoming edges onto outgoing ones and vice versa, so the numbers of incoming and outgoing
edges would be equal, which is impossible for degree $3$.
\item The centre is an edge $e=(u\to v)$. Then $\sigma(e)=e$ and $\sigma$ reverses it,
so $\sigma$ interchanges $u$ and $v$. Hence the branch at $v$ is the reversal of the branch at $u$.
Conversely, for any branch $X$ of type $-$ the tree obtained by gluing $X$ to $\overline X$ along an edge
has an anti-automorphism. Such trees correspond bijectively to branches $X$ of type $-$
with $n/2$ saddles.
\end{itemize}
Therefore $F_n=K_{n/2}$ for even $n$ and $F_n=0$ for odd $n$.
\end{proof}

\begin{table}[ht]
\centering\small
\begin{tabular}{c|ccccccccccccc}
$n$ & 0 & 1 & 2 & 3 & 4 & 5 & 6 & 7 & 8 & 9 & 10 & 11 & 12\\\hline
$K_n$ & 1 & 2 & 6 & 25 & 111 & 540 & 2736 & 14396 & 77649 & 427608 & & & \\
$T_n$ & 1 & 2 & 4 & 14 & 48 & 198 & 853 & 3952 & 18971 & 94428 & 481898 & 2513964 & 13347066\\
$H_n$ & 1 & 1 & 3 & 7 & 27 & 99 & 439 & 1976 & 9541 & 47214 & 241219 & 1256982 & 6674901
\end{tabular}
\caption{Numbers of Reeb trees with $n$ saddles. The values of $T_n$ and $H_n$
for $n\le 8$ agree with a direct computer enumeration of all trees.}
\label{tab:trees}
\end{table}

\begin{remark}
Equation \eqref{eq:K} is of Pólya--Otter type, so one expects a square-root
singularity and $T_n\sim C\rho^{\,n}n^{-5/2}$. We only observed this numerically:
$T_n/T_{n-1}\approx 5.98$ for $n=30$. Since $F_n$ grows like $\rho^{n/2}$, this
suggests $H_n\sim T_n/2$.
\end{remark}

\section{Topology of Functions with Isolated Critical Points on a Surface}

\subsection{Local Classification}

Let $F$ be a closed smooth surface, and $f: F \rightarrow {\mathbb R}$ be a smooth function with a finite number of critical points. Due to the compactness of the surface, this condition is equivalent to each critical point being isolated.

\begin{lemma} If $y$ is a regular value of the function $f$, then $f^{-1}(y)$ consists of a finite number of circles. For any open set $U$ with a boundary that is transverse to $f^{-1}(y)$, the intersection $f^{-1}(y) \cap U$ has a finite number of connected components.
\end{lemma}

\begin{proof} By contradiction. Suppose $f^{-1}(y)$ contains an infinite number of circles. Then they partition the surface into an infinite number of regions, each containing a local minimum or maximum, which contradicts the finiteness of the number of critical points.

If $f^{-1}(y) \cap U$ has an infinite number of connected components, then there are also infinitely many intersection points $f^{-1}(y) \cap \partial U$. Thus, there exists a point $x \in f^{-1}(y) \cap \partial U$ that is a boundary point of the set $f^{-1}(y) \cap \partial U$. The derivatives of the function $f$ at the point $x$ along $f^{-1}(y)$ and $\partial U$ are both equal to 0. Therefore, the point $x$ is critical. This contradicts the conditions of the lemma.
\end{proof}

\begin{lemma} {For each critical point $x_{0}$ with critical value} $y_{0} = f(x_{0})$, {there exists a neighborhood} $U(x_{0})$, {in which}
$$f^{-1}(y) \cap {\text{cl}}(U(x_{0})) \cong {\text{Con}}(\cup S^{0}).$$
Here, Con $(\cup S^{0})$ is the cone over the union of 0-dimensional spheres, that is, the union of an even number of segments, which share a common point—the endpoint of each segment.
\end{lemma}

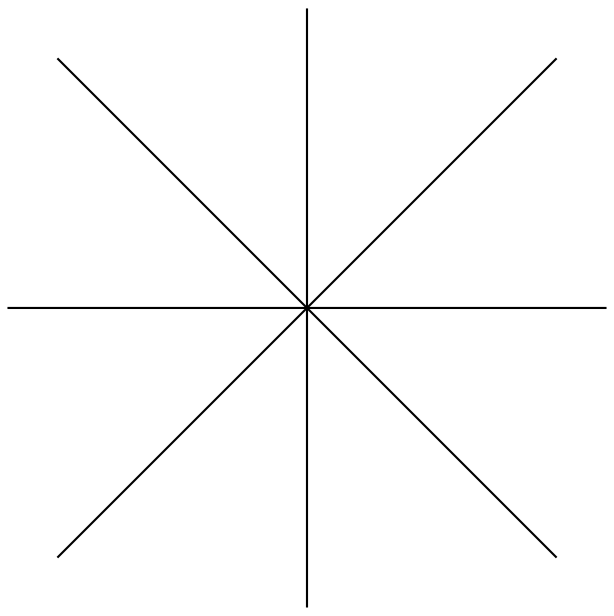
\begin{figure}[ht]

\center{ \begin{tikzpicture}[x=0.1\textwidth, y=0.1\textwidth]
  \draw (-1.2,0) -- (1.2,0);      
  \draw (0,-1.2) -- (0,1.2);      
  \draw (-1,-1) -- (1,1);     
  \draw (-1,1) -- (1,-1);     
\end{tikzpicture}}
\caption{Critical level in the neighborhood of a critical point}
\label{p11}
\end{figure}

\begin{proof} Since $x_{0}$ is an isolated critical point, there exists a neighborhood $U$ that does not contain other critical points of the function $f$. Let $g$ be the restriction of the function $f$ to the set $U \setminus \{x_{0}\}$. Then, on this set, the level curves of both functions coincide. Since $y_{0}$ is a regular value of the function $G$, the preimage $g^{-1}(y_{0})$ is a one-dimensional submanifold in $U \setminus \{x_{0}\}$, that is, a disconnected union of closed curves and open intervals.

Without loss of generality, we can assume that $\partial U$ intersects transversally with $f^{-1}(y_{0})$ and does not contain critical points. Similarly to Lemma 14.1, it can be shown that the intersection points are finite in number.

We  restrict the neighborhood $U$ to a neighborhood $U(x_{0})$, the closure of which does not intersect with the circles and arcs where both ends lie on $\partial U$, and such that $\partial U(x_{0})$ intersects transversally with $g^{-1}(y_{0})$ and does not contain critical points. Then, $f^{-1}(y_{0}) \cap cl(U(x_{0}))$ consists of a union of closed arcs, at least one of whose ends coincides with the point $x_{0}$.

We  show that none of these arcs form a loop (i.e., both ends do not coincide with the point $x_{0}$). By contradiction, suppose such an arc exists; then it bounds a two-dimensional disk, and the function $f$ takes a minimum or maximum value on it that is different from $y_{0}$. Thus, there exists another critical point (or minimum or maximum) in $U(x_{0})$ that is different from the point $x_{0}$. This contradicts the choice of the neighborhood $U(x_{0})$.

\begin{theorem} {For every critical point $x _{0}$ (except local minima and maxima) there exists a neighborhood in which the function is conjugate to the function} Re $z^{k}$ for some number $k \in {\mathbb Z_+}$.
\end{theorem}
The level curves of the function Re $z^{k}$ for $k = 3$ are schematically depicted in Figure \ref{p12}.

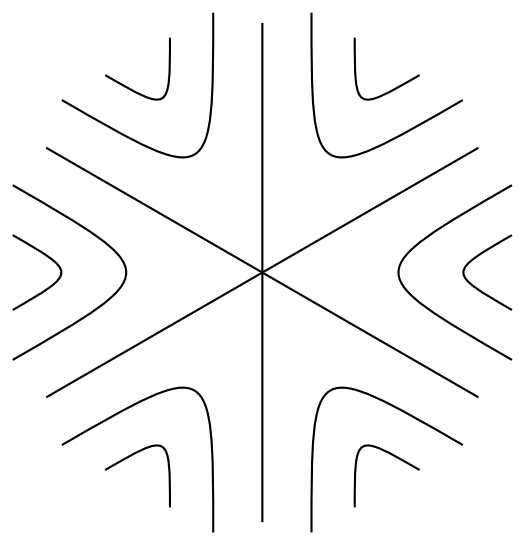
\begin{figure}[ht]
\center{\begin{tikzpicture}[x=0.1\textwidth, y=0.1\textwidth]
  \foreach \a in {30,90,150}
    \draw (\a+180:1) -- (\a:1);

  \foreach \k in {0,...,5} {
    \begin{scope}[rotate=60*\k]
      \foreach \q in {0.15,0.35} {
        \draw (1,-\q) .. controls ({1-1.732*\q},0) .. (1,\q);
      }
    \end{scope}
  }
\end{tikzpicture} }
\caption{The stratification of the neighborhood of the critical level}
\label{p12}
\end{figure}

We  show that there is a finite number of arcs. Suppose, for contradiction, that there are infinitely many. We cut the surface $f$ along $f^{-1}(y _{0})$. Since cutting along each closed curve with endpoints on $f^{-1}(y _{0})$ decreases the genus of the surface or increases the number of connected components, we obtain an infinite number of connected components. On each connected component, there exists a point of local minimum or maximum. This contradicts the finiteness of the number of critical points on the surface $f$. The resulting contradiction proves that the critical level in the vicinity of the critical point is a union of a finite number of arcs with a common endpoint.

We  show that the number of arcs is even. We  call the components of the set $ U(x_{0})\setminus f^{-1}(y _{0})$ with $f(x) > y _{0}$ positive regions, and the components with $f(x) < y _{0}$ negative regions. Each arc lies on the boundary of both a positive and a negative region. These regions alternate. Thus, there are an equal number of regions of each type, and the total number of regions, as well as arcs, is even.
\end{proof}

\begin{proof} See \cite{prishlyak2002topological1}. 
\end{proof}

\noindent \textbf{Corollary 1.} If $k=1$, then there exists a neighborhood of the critical point in which the function is conjugate to a function without critical points.

\noindent \textbf{Remark.} Similarly to the proof of the theorem, it can be shown that for every local minimum (maximum), there exists a neighborhood where the function is conjugate to the function
$$g(x,y)=x ^{2}+y ^{2} (g(x,y)= -x ^{2}-y^{2}).$$

\begin{definition} We call the Poincaré index ind$^{f}(x _{0})$ of the critical point $x _{0}$ of the function $f$ as the Poincaré index of its gradient field grad $f$ for some Riemannian metric.
\end{definition}

\noindent \textbf{Consequence 2.} Let $x _{0}$ be an isolated critical point of the function $f$, and $y _{0}$ be a critical point of the function $g$. The functions $f$ and $g$ are topologically equivalent in some neighborhoods of these points that do not contain other critical points if and only if they have the same Poincaré indices, ind$^{f}(x _{0}) =$ ind$^{g}(y _{0})$.

\noindent \textbf{Remark.} For critical points on a surface, ind$^{f}(x _{0}) \leq 1$. For the function $f = \text{Re} \, z^{k}$, the Poincaré index ind$^{f}(0) = 1 - k$.

\subsection{ Global Classification}

Let $F$ be a smooth closed surface and $f: F \rightarrow {\mathbb R}$ be a smooth function with isolated critical points and critical values $y _{1}, \dots , y _{n}$ such that $y_{i} < y_{j}$ if $i < j$. By Theorem 1, the preimage of each critical value is homeomorphic to a graph in union with closed curves. The edges of the graph and the closed curves are smoothly embedded 1-manifolds. On each closed curve, we fix a point, thereby transforming this curve into a loop with a distinguished vertex.

Let $G _{i}(f)$ denote the graph (possibly disconnected) that geometrically coincides with the preimage of the critical value $y _{i}$ of the function $f$. The vertices of the graph are the critical points and the distinguished points.

\begin{definition} The surface $F$ together with the embedded graphs $G _{i}(f)$ is called the {diagram} $D = \{F, G _{1}(f), \dots , G_{n}(f)\}$ of the function $f$.
\end{definition}

\begin{definition} Two diagrams are called {isomorphic} if there exists a homeomorphism of the surfaces that maps graphs to graphs, and vertices of the graphs to vertices, preserving the order of the graphs.
\end{definition}

\begin{theorem} {Two functions with isolated critical points on closed surfaces are topologically equivalent if and only if their diagrams are isomorphic.}
\end{theorem}

\begin{proof} Necessity. If there is a given homeomorphism of the surface that connects, then it maps critical levels to critical levels and thus defines an isomorphism of diagrams.

Sufficiency. Let the diagrams of two functions be isomorphic. Then there exists a homeomorphism $g$ of the surface $F$ that maps critical levels to critical levels. We cut the surface $f$ along the critical levels. We obtain surfaces $F_{k}$ that are homeomorphic to cylinders $S^{1} \times [0,1]$. Then the homeomorphism $g$ induces homeomorphisms of the obtained cylinders. We replace these homeomorphisms with homeomorphisms that map the levels of the function to levels proportionally to the value of the function between two critical values and coincide with the original on each of the bases of the cylinders. Since the constructed homeomorphisms coincide on the boundaries, they define the sought homeomorphism of the surface.
\end{proof}
Let $D =\{F, G _{1}(f), \dots , {G_{n}}(f)\}$ be the diagram of the function $f$ with isolated critical points on the surface $F$. Each cylinder $F_{k}$ with bases on the graphs $G _{i}(f), G_{i+1}(f)$ defines one cycle corresponding to the base of the cylinder on these graphs. The cycle on $G_{i}(f)$ is called the upper cycle, and on $G _{i+1}(f)$ --- the lower cycle. In this case, local minima and maxima correspond to cycles consisting of a single vertex. It is worth noting that each edge is contained in only one upper and only one lower cycle.

\begin{definition} A \textit{distinguishing graph} is a graph $G$ partitioned into subgraphs $G_{1},\dots , G _{n}$, on which cycles are highlighted that are divided into two groups --- upper and lower, and a one-to-one correspondence $\psi$ is defined between the lower cycles from $G _{i}$ and the upper cycles from $G _{i+1}$, $i=1,\dots ,n-1$. In this case, each edge is contained in only one upper and one lower cycle, and each isolated vertex forms one upper or lower cycle.

An isomorphism of distinguishing graphs is an isomorphism of graphs such that subgraphs are mapped to subgraphs, upper cycles to upper cycles, lower cycles to lower cycles, and there is also a correspondence preserved between upper and lower cycles.
\end{definition}
Every function with isolated special points on a closed surface defines a unique distinguishing graph, up to isomorphism.

\begin{theorem} {Two functions are topologically conjugate if and only if their distinguishing graphs are isomorphic.}
\end{theorem}

\begin{proof} Necessity. If the functions are topologically conjugate, then their diagrams are isomorphic. The restriction of this isomorphism to the graph defines an isomorphism of distinguishing graphs.

Sufficiency. Suppose the distinguishing graphs are isomorphic. We  construct an isomorphism of diagrams. By construction, each pair of upper and lower cycles, related by a one-to-one correspondence $\psi$, corresponds to some cylinder $F_{k}$. Thus, a correspondence between cylinders is established. We  construct homeomorphisms between the cylinders in such a way that their restrictions to the upper and lower bases coincide with the restrictions of the graph isomorphisms to the corresponding cycles. Since these homeomorphisms  coincide on the boundaries, they define a homeomorphism of the surface $F$ and, consequently, an isomorphism of diagrams.
\end{proof}

Now we investigate the question of when a distinguishing graph defines a function with isolated special points on a closed surface. Two edges $e_{1}, e_{2}$, incident to a common vertex $v$, are called adjacent for the vertex $v$ if there exists a cycle in the distinguishing graph that contains the fragment $(e_{1}, v, e_{2})$ or $(e_{2}, v, e_{1})$.

\begin{theorem} {A distinguishing graph defines a function with isolated special points on a closed surface if and only if for any two edges incident to a common vertex, there exists a sequence of edges in which each subsequent edge is adjacent to the previous one at that vertex.}
\end{theorem}

\begin{proof} The necessity directly follows from the construction. We  prove the sufficiency. The conditions of the theorem guarantee that after attaching cylinders to the graph, we  obtain a closed surface. We  construct a function $f$ on it. For this, we set $f(G_{i}) = i$, and for the cylinder $S^{1} \times [0,1]$ with bases on $G_{i}$ and $G_{i+1}$, we have $$f(S^{1} \times \{ t \}) = i + t.$$ In this way, at the vertices of the graph, the cylinders are attached such that the function $f$ coincides with the function $i + \text{Re} \ z^{k}$ in the corresponding coordinate system.
\end{proof}

\begin{remark}
The cell decompositions and rotation systems used in this section can be
enumerated and visualized with~\cite{prishlyak2026rotationcodes,
prishlyak2026cellembedded}; chord diagrams with one cycle, which describe
one-vertex embeddings, with~\cite{prishlyak2026chord}; the corresponding handle
decompositions with~\cite{prishlyak2026handle}.
\end{remark}

\subsection{Functions with three critical points}\label{sec:three-crit}

For a function with three critical points, one of them is a local minimum,
another is a local maximum, and the third one is the vertex of a bouquet of
circles, which is a critical level of the function. According to the previous
sections, the topological structure of such a function is determined by the
following cell decomposition of the surface: the $1$-skeleton is a bouquet
$B_k$ of $k$ circles with the vertex $v$, one $2$-cell contains the maximum
point, and the other $2$-cell contains the minimum point. Below we describe all
such cell decompositions by means of signed permutations.

\paragraph{Encoding.}
Both $2$-cells are $k$-gons whose sides are glued to the circles of $B_k$ by
linear homeomorphisms, and every circle of $B_k$ lies on the boundary of both
$2$-cells (each circle is covered exactly once by each polygon). Number the
sides of the first $k$-gon cyclically along its boundary by $1,2,\dots,k$.
This numbers the circles of $B_k$ and orients them (the orientation of the
$i$-th circle is the one induced from the side $i$ of the first polygon).
Now traverse the boundary of the second $k$-gon, starting from the side that
is glued to the circle $1$, in the direction in which this circle is passed
along its orientation. Writing down the numbers of the circles in the order
of traversal, and attaching the sign $\varepsilon=+$ if the circle is passed
along its orientation and $\varepsilon=-$ otherwise, we obtain a signed
permutation
\[
  \pi=\bigl(\pi_1^{\varepsilon_1},\pi_2^{\varepsilon_2},\dots,\pi_k^{\varepsilon_k}\bigr),
  \qquad \pi_1=1,\ \varepsilon_1=+ .
\]
There are $(k-1)!\,2^{k-1}$ such words, and the word $\pi$ determines the
cell complex $X_\pi$ uniquely: it is obtained from $B_k$ by attaching the
first $k$-gon along $1,2,\dots,k$ and the second one along $\pi$.

\paragraph{Local planarity at the vertex.}
The complex $X_\pi$ is a closed surface if and only if the vertex $v$ has a
neighborhood homeomorphic to a disc, i.e.\ if the link of $v$ is a single
circle (at the interior points of the edges and the $2$-cells the surface
condition holds automatically). The link $L_\pi$ is the graph defined as
follows. Its $2k$ vertices are the ends $i^-,i^+$ (tail and head) of the
$k$ circles. Its $2k$ edges are the corners of the two polygons:
\begin{itemize}
  \item the first polygon gives the corners joining $i^+$ with $(i+1)^-$,
        $i=1,\dots,k$ (indices modulo $k$);
  \item for two consecutive entries $\pi_j^{\varepsilon_j},\pi_{j+1}^{\varepsilon_{j+1}}$
        of the second word (cyclically), the corner joins the end of the circle
        $\pi_j$ at which we leave it (the head $\pi_j^+$ if $\varepsilon_j=+$,
        the tail $\pi_j^-$ if $\varepsilon_j=-$) with the end of the circle
        $\pi_{j+1}$ at which we enter it (the tail if $\varepsilon_{j+1}=+$,
        the head if $\varepsilon_{j+1}=-$).
\end{itemize}
Every end of every circle belongs to exactly two corners, one in each polygon,
so $L_\pi$ is a $2$-regular graph, i.e.\ a disjoint union of cycles. Hence
$X_\pi$ is a surface if and only if $L_\pi$ is connected. This is checked in
linear time with a union--find structure.

\paragraph{The resulting surface.}
If $L_\pi$ is connected, $X_\pi$ is a closed connected surface with one vertex,
$k$ edges and two faces, so
\[
  \chi(X_\pi)=1-k+2=3-k .
\]
The surface is orientable if and only if all signs are $+$. Indeed, an
orientation of $X_\pi$ is a choice of orientations of the two $k$-gons such
that every circle is traversed in opposite directions by the induced boundary
orientations. The circle $i$ is traversed by the first polygon along its
orientation, and by the second one with the sign $\varepsilon$, so a
compatible choice exists if and only if all the signs of $\pi$ are equal,
that is, all $\varepsilon_i=+$ (the first polygon is then oriented
counterclockwise and the second one clockwise). Consequently an orientable
surface occurs only for odd $k=2g+1$ and has genus $g=(k-1)/2$, while a
non-orientable one has $k-1$ crosscaps. For orientable surfaces the number
of admissible words is the number of $k$-cycles $\beta$ such that $\alpha\beta$
is again a $k$-cycle for $\alpha=(1\,2\dots k)$; by a classical result it equals
$2(k-1)!/(k+1)$ for odd $k$ (and $0$ for even $k$). We have used it as a test of
the program below.

\paragraph{Equivalence.}
The numbering of the first polygon depends on the choice of the first side
and of the direction of traversal. Changing them by an element of the dihedral
group $D_k$ gives another word that describes the same cell decomposition.
Explicitly, for the rotation by $r$ the circle $i$ gets the number $i-r \pmod k$
and the signs stay the same; for the reflection the circle $i$ gets the number
$-i-r \pmod k$ and all the orientations of the circles are reversed, so all
signs change. After such a renumbering the second word is cyclically shifted
so that it starts with the circle $1$; if the sign of this entry is $-$, the
boundary of the second polygon is read in the opposite direction (the order of
the entries is reversed and all signs are changed), which makes the first sign
equal to $+$ again. Two words are \emph{equivalent} if they are in the same
$D_k$-orbit. Equivalence classes correspond to the topologically different
cell decompositions with the distinguished maximum and minimum cells.

\paragraph{Algorithm.}
The program below enumerates all words $\pi$ with $\pi_1=1$,
$\varepsilon_1=+$, keeps those for which $L_\pi$ is connected, and splits them
into orientable and non-orientable ones. For every kept word that has not been
met yet, its whole $D_k$-orbit (normalized as above) is put into a set of
visited words, and the counter of classes is increased by one. The running time
is $O\bigl((k-1)!\,2^{k-1}\cdot k\bigr)$, which is feasible up to $k\approx 9$.

\begin{lstlisting}[language=Python,basicstyle=\ttfamily\footnotesize,breaklines=true]
import sys
from itertools import permutations, product

def link_is_circle(word, k):
    """word = [(edge, sign), ...], edges 0..k-1, sign = +1/-1.
    Returns True iff the link of the vertex is connected (one cycle)."""
    parent = list(range(2 * k))      # node 2i = tail of loop i, 2i+1 = head
    def find(x):
        while parent[x] != x:
            parent[x] = parent[parent[x]]
            x = parent[x]
        return x
    def union(a, b):
        parent[find(a)] = find(b)
    for i in range(k):               # corners of the first polygon
        union(2 * i + 1, 2 * ((i + 1) % k))
    for j in range(k):               # corners of the second polygon
        (e, s), (e2, s2) = word[j], word[(j + 1) % k]
        leave = 2 * e + (1 if s > 0 else 0)    # end of loop e where we leave
        enter = 2 * e2 + (0 if s2 > 0 else 1)  # end of loop e2 where we enter
        union(leave, enter)
    return len({find(x) for x in range(2 * k)}) == 1

def normalize(word):
    """Start at edge 0; if it carries sign -, read the boundary backwards."""
    i = [e for e, _ in word].index(0)
    w = word[i:] + word[:i]
    if w[0][1] < 0:
        w = [(0, 1)] + [(e, -s) for e, s in reversed(w[1:])]
    return tuple(w)

def dihedral_orbit(word, k):
    """All words obtained by renumbering the first polygon (D_k action)."""
    for r in range(k):
        yield normalize([((e - r) % k, s) for e, s in word])      # rotation
        yield normalize([((-e - r) % k, -s) for e, s in word])    # reflection

def count(k):
    """Returns (orientable words, orientable classes,
                nonorientable words, nonorientable classes)."""
    seen = set()
    ow = oc = nw = nc = 0
    for perm in permutations(range(1, k)):
        order = (0,) + perm
        for signs in product((1, -1), repeat=k - 1):
            word = tuple(zip(order, (1,) + signs))
            if not link_is_circle(word, k):
                continue
            orientable = all(s > 0 for _, s in word)
            if orientable: ow += 1
            else:          nw += 1
            if word not in seen:
                seen.update(dihedral_orbit(list(word), k))
                if orientable: oc += 1
                else:          nc += 1
    return ow, oc, nw, nc

if __name__ == "__main__":
    for k in range(1, int(sys.argv[1]) + 1):
        print(k, 3 - k, *count(k), flush=True)
\end{lstlisting}

\paragraph{Results.}
The numbers of admissible words and of equivalence classes are given in
Table~\ref{tab:three-crit}. For orientable surfaces the numbers of classes
$1,1,4,25$ ($k=1,3,5,7$) and for non-orientable ones $1,1,4$ ($k=2,3,4$)
agree with the direct count.

\begin{table}[ht]
\centering
\begin{tabular}{rrrrrr}
\toprule
$k$ & $\chi$ & \multicolumn{2}{c}{orientable} & \multicolumn{2}{c}{non-orientable}\\
    &        & words & classes & words & classes\\
\midrule
1 &  $2$ &   1 &  1 &      0 &     0\\
2 &  $1$ &   0 &  0 &      1 &     1\\
3 &  $0$ &   1 &  1 &      3 &     1\\
4 & $-1$ &   0 &  0 &     20 &     4\\
5 & $-2$ &   8 &  4 &    140 &    20\\
6 & $-3$ &   0 &  0 &   1348 &   125\\
7 & $-4$ & 180 & 25 &  14924 &  1136\\
8 & $-5$ &   0 &  0 & 198144 & 12504\\
\bottomrule
\end{tabular}
\caption{Cell decompositions of closed surfaces with one vertex, $k$ edges and two
faces: admissible signed permutations and their $D_k$-classes.}
\label{tab:three-crit}
\end{table}


\section{Literature Review}

Classical textbooks on dynamical systems, structural stability and
Morse--Smale systems include \cite{Katok1995, Meiss2017, Palis1982}.
Morse theory is presented in \cite{Milnor1963, prish2002Morse}. Simplicial
complexes, Euler characteristics and homology groups are fundamental concepts of
homology theory \cite{hat02, prish-modtop06}. Regular cell complexes and
discrete Morse functions are discussed in discrete Morse theory \cite{for98}.
SLW-graphs were introduced in \cite{prishlyak1997graphs}. Rotation systems are
one of the core concepts of topological graph theory
\cite{Peixoto1973, Gross1987, mohar2001graphs}, and the enumeration of chord
diagrams is covered in \cite{Kadubovski2021}.

The Reeb graph, the space of connected components of the level sets of a
function, was introduced by G.~Reeb~\cite{Reeb1946} for Morse functions, and
A.\,S.~Kronrod~\cite{Kronrod1950} studied the structure of the level sets of
continuous functions of two variables, which for functions on the sphere leads
to trees (the Kronrod--Reeb tree). These works contain the classification of
simple Morse functions by Reeb graphs and the realization theorem for graphs, which we
use in Section~\ref{sec:reeb}. A.\,T.~Fomenko and A.\,V.~Bolsinov developed
the invariants of integrable Hamiltonian systems with one degree of freedom and
of Morse functions on surfaces in terms of atoms and molecules, that is, Reeb
graphs with additional marking of the vertices~\cite{Bolsinov2004}.
V.\,V.~Sharko obtained the global topological classification of Morse functions
on surfaces and studied the algebraic and topological aspects of functions on
manifolds~\cite{Sharko1993}. 


The structures considered in this series are accompanied by open software.
Morse flows on the $2$-disc with holes can be specified by codes and
enumerated and visualized with the program~\cite{prishlyak2026morsedisk};
handle decompositions of surfaces given by rotation systems are drawn
in~\cite{prishlyak2026handle}, and Morse--Smale (MS), SN and SC flows and their bifurcations
are visualized from rotation-system codes in~\cite{prishlyak2026mssnsc}.
Rotation systems of cell-embedded graphs are enumerated, encoded and
visualized in~\cite{prishlyak2026rotationcodes, prishlyak2026cellembedded},
and one-cycle signed chord diagrams in~\cite{prishlyak2026chord}.

Topological properties of functions and vector fields have been studied in many works. Among them, we would like to highlight the works of the Kyiv topological school: 
textbooks \cite{prish2002Morse, prish-modtop06, prish2015top, prish23algkom}, 
scientific articles by the author on function topology \cite{Prishlyak1993, prishlyak1998, prishlyak1999equivalence, prishlyak2000conjugacy, prishlyak2001conjugacy, prish2001top, prishlyak2002topological1, prishlyak2002morse, prishlyak2003regular, prishter2024}, vector fields \cite{Prish1997vec, Prishlyak2001, prishlyak2002morse1, Prishlyak2002, prishlyak2003topological, prishlyak2003sum, prishlyak2005complete, Prishlyak2007}, and other geometric objects \cite{Prishlyak1994, prishlyak1997graphs, Prishlyak1999}, as well as the works of his students: K. Myshchenko \cite{prishlyak2007classification}; N. Budnytska \cite{Bud2008knu}; D. Lychak \cite{lychak2009morse}; A. Bondarenko \cite{Bond2012mfat}; O. Vyatychaninova \cite{VyatP2013Mol}; Bohdana Hladysh \cite{Hladysh2016, hladysh2017topology,
Hladysh2019, PLH2023}; A. Prus \cite{Prishlyak2017, prishlyak2020three, Prishlyak2021}; V. Kiosak \cite{KPL2022}; S. Bilun \cite{bilun2002closed, PBP23}
V. Lisikevich \cite{LisP2013KNU}, I. Ivaniuk \cite{IvanPrish2014-func-deform3, IP2015knu},\cite{Ovtsynov2024}, D. Skotchko \cite{
PS2016-F-atoms}, M. Loseva \cite{Losieva2017, Prishlyak2019, Prishlyak2020}, Z. Kybalko \cite{Kybalko2018}, K. Khatamian \cite{Hatamian2020}.


\section*{Software availability}
The programs related to this paper are archived on Zenodo with permanent DOIs:
the viewer and enumerator of Morse flow codes on the $2$-disc with
holes~\cite{prishlyak2026morsedisk}, the Smale handle decomposition viewer for
surfaces (also available online at \url{https://prishlyak.github.io/smale-handle-viewer/})~\cite{prishlyak2026handle},
the visualizer of MS, SN and SC flows and their bifurcations~\cite{prishlyak2026mssnsc},
the enumeration and visualization of rotation-system codes of cellular graph
embeddings~\cite{prishlyak2026rotationcodes, prishlyak2026cellembedded} and of one-cycle signed chord
diagrams~\cite{prishlyak2026chord}.

\section{Conclusion}

We studied the topological equivalence of smooth functions and of Hamiltonian
flows on surfaces and the discrete invariants that classify them. The main results are
the following.

\begin{enumerate}
\item We distinguished two equivalence relations for functions: the strong
equivalence $f\circ h=h'\circ g$ with an increasing $h'$, which keeps track of the order of all
critical values, and the level equivalence, which only requires the components of
the level sets to be mapped onto each other with the direction of increase preserved.
Simple Morse functions on closed orientable surfaces are classified by Reeb graphs
with the order of vertices up to the first relation and by oriented Reeb graphs up to
the second; the second is the equivalence that corresponds to Hamiltonian
flows. Simple Hamiltonian flows are classified by oriented Reeb graphs up to isomorphism and the simultaneous reversal
of all orientations, and for trees by trees with one marked edge at every vertex of degree three.
We also gave a complete classification of functions with three and four
extrema on the circle and on the line.

\item For the disc we obtained the recurrence for the number $K_n$ of rooted
Reeb trees with $n$ saddles, equivalently the equation
$K(x)=1+\tfrac{x}{2}\bigl(3K(x)^2+K(x^2)\bigr)$ for its generating function.

\item For the sphere the generating function of the numbers $T_n$ of simple
Morse functions with $n$ saddles up to level equivalence is
$K(x)^2\bigl(1-2xK(x)\bigr)$, and the number of Hamiltonian flows is
$H_n=\tfrac12\bigl(T_n+F_n\bigr)$, where $F_n=K_{n/2}$ for even $n$ and $F_n=0$
for odd $n$. The first values for $n=0,\dots,7$ are
\[
T_n=1,\,2,\,4,\,14,\,48,\,198,\,853,\,3952,\qquad
H_n=1,\,1,\,3,\,7,\,27,\,99,\,439,\,1976 .
\]
Numerically the ratio $T_n/T_{n-1}$ approaches approximately $6$ (it equals
$5.98$ for $n=30$), which is consistent with the square-root singularity expected for
equations of this type; we did not prove the asymptotics.

\item For functions with isolated, possibly degenerate, critical points we described
the topological structure by cell decompositions of the surface. For functions with three
critical points (a minimum, a maximum and a critical point whose level set is a bouquet
of $k$ circles) we gave an algorithm based on signed permutations and
computed the numbers of classes with respect to the dihedral group $D_k$:
$1,1,4,25$ for $k=1,3,5,7$ on the orientable surfaces (genus $0,1,2,3$) and
$1,1,4$ for $k=2,3,4$ on the non-orientable ones.
\end{enumerate}

\medskip
\noindent\textbf{Further directions.}
Several questions remain open. The first one is the exact asymptotics of $T_n$ and $H_n$
(the radius of convergence of $K(x)$ and the constant in $T_n\sim C\rho^{\,n}n^{-5/2}$).
The second is the enumeration of the classes of simple Morse functions up to strong
equivalence, i.e.\ of Reeb trees with the order of the critical values, which
is finer than the counts given here. The third is the enumeration of
Reeb graphs with cycles, for orientable surfaces of positive genus, by the same
combination of rooted trees and Burnside's lemma, and of functions on
non-orientable surfaces, where the Reeb graph is not a complete invariant.
Finally, it would be natural to extend the computations for functions with
isolated critical points beyond three critical points and to compare the numbers
obtained with those arising from the enumeration problems for Morse functions
considered earlier in the literature.


\section*{Declaration of generative AI use}

During the preparation of this work the author used the large language model
Claude (Anthropic, Claude Sonnet 5.5, via the claude.ai web interface,
November 2026) as a tool for the following tasks: drafting and editing the
English text of the abstract, the introduction, the conclusion and parts of
Sections on topological equivalence,  Reeb graphs and Hamiltonian flows;  and writing auxiliary Python programs that enumerate trees and check the formulas by direct enumeration for small numbers of saddles. The numerical results given in the paper were obtained by programs, and the formulas were compared with independent enumeration. The author reviewed, corrected and verified all the AI-assisted text, code and references, and takes full responsibility for the content of
the publication. The AI tool is not an author of this work. The software cited in this paper was written by the author.

\bibliographystyle{plainnat}
\bibliography{prish}
\end{document}